\documentclass[12pt]{amsart}

\usepackage{amsmath,amssymb,amsthm,amsfonts}

\newtheorem{thm}{Theorem}
\newtheorem{prop}{Proposition}
\newtheorem{defin}{Definition}

\newtheorem{cor}{Corollary}
\newtheorem{lem}{Lemma}
\newtheorem{conj}{Conjecture}

\newcommand{\CC}{\mathfrak C}

\newcommand{\G}{{\mathfrak G}}

\newcommand{\hh}{{\mathfrak H}}

\newcommand{\N}{{\mathbb N}}

\newcommand{\Q}{{\mathbb Q}}

\newcommand{\R}{{\mathbb R}}

\newcommand{\TT}{{\mathbb T}}

\newcommand{\V}{{\mathfrak V}}

\newcommand{\Z}{{\mathbb Z}}

\newcommand{\blo}{\text{block}} 

\newcommand{\id}{\text{Id}}   

\newcommand{\Log}{\text{Log}}

\newcommand{\mac}{\text{SL}(2,\R)}   
\newcommand{\mcz}{\text{SL}(2,\Z)}   
\newcommand{\macn}{\text{SL}(n,\R)}   
\newcommand{\mczn}{\text{SL}(n,\Z)}   

\newcommand{\Sqrt}{\text{Sqrt}}     

\newcommand{\ga}{\gamma}\newcommand{\Ga}{\Gamma}
\newcommand{\de}{\delta}

\newcommand{\Om}{\Omega}

\newcommand{\tB}{\tilde{B}}

\newcommand{\tm}{\tilde{m}}\newcommand{\tM}{\tilde{M}}

\newcommand{\va}{\vec{a}}
\newcommand{\vb}{\vec{b}}
\newcommand{\vde}{\vec{\de}}
\newcommand{\vk}{\vec{k}}
\newcommand{\vl}{\vec{l}}
\newcommand{\vvp}{\vec{p}}
\newcommand{\vq}{\vec{q}}

\newcommand{\vx}{\vec{x}}
\newcommand{\vy}{\vec{y}}

\newcommand{\tGa}{\tilde{\Ga}}

\date{\today}

\begin{document}

\title[Blocking in homogeneous spaces]
{Connection blocking in nilmanifolds and other homogeneous spaces}

\author{Eugene Gutkin}

\address{Nicolaus Copernicus University, Department of Mathematics, Chopina 12/18,
87-100 Torun and Mathematics Institute of the Polish Academy of
Sciences, Sniadeckich 8, 00-956 Warsaw, Poland}

\email{gutkin@mat.umk.pl,\ gutkin@impan.pl}


\date{\today}

\begin{abstract}
Let $G$ be a connected Lie group acting locally simply
transitively on a manifold $M$. By {\em connecting curves} in $M$
we mean the orbits of one-parameter subgroups of $G$. To {\em
block} a pair of points $m_1,m_2\in M$ is to find a {\em finite
set} $B\subset M\setminus\{m_1,m_2\}$ such that every connecting
curve joining $m_1$ and $m_2$ intersects $B$. The homogeneous
space $M$ is {\em blockable} if every pair of points in $M$ can be
blocked. Motivated by the {\em geodesic security} \cite{BanGut10},
we conjecture that the only blockable homogeneous spaces of finite
volume are the tori $\R^n/\Z^n$. Here we establish the conjecture
for {\em nilmanifolds}.
\end{abstract}

\maketitle

\tableofcontents

\section{Introduction}  \label{intro}
The theme of finite blocking has its genesis in a Leningrad
Mathematical Olympiad problem \cite{FoKi94,Wi07} worded as
follows. The president and a terrorist are moving in a rectangular
room. The terrorist intends to shoot the president with his `magic
gun' whose bullets bounce of the walls perfectly elastically: The
angles of incidence and reflection are equal. Presidential
protection detail consists of superhuman body guards. They are not
allowed to be where the president or the terrorist are located,
but they can be anywhere else, changing their locations
instantaneously, as the president and the terrorist are moving
about the room. Their task is to put themselves in the way of
terrorist's bullets shielding the president. The problem asks how
many body guards suffice.

To translate this into mathematical setting, let $\Om$ be a
bounded plane domain. For arbitrary points $p,t\in\Om$ let
$\Ga(p,t)$ be the family  of billiard orbits in $\Om$ connecting
these points. Body guards correspond to
$b_1,\dots,b_N\in\Om\setminus\{p,t\}$ such that every
$\ga\in\Ga(p,t)$ passes through one of these points. If for any
$p,t\in\Om$ there is a {\em blocking set}
$B=B(p,t)=\{b_1,\dots,b_N\}$ then the domain is {\em uniformly
secure}. The minimal possible $N$ is then the {\em blocking
number} of $\Om$. The Olympiad problem is to show that a rectangle
is uniformly secure and to find its blocking number.

The solution is an exercise in plane geometry based on two facts:
i) A rectangle tiles the euclidean plane under reflections; ii)
The torus $\TT^2=\R^2/\Z^2$ is uniformly secure, where the role of
billiard orbits is played by the images of straight lines under
the projection $\R^2\to\TT^2$. A blocking set in the torus is the
set of midpoints of all joining segments: It comprises at most $4$
points. A blocking set in the rectangle is also the set of
midpoints of all joining billiard orbits. It comprises at most
$16=4\times 4$ points where the factor $4$ is due to the $4$
copies of the rectangle needed to tile the torus.

The billiard orbits in the rectangle and the straight lines in the
torus are examples of geodesics in riemannian manifolds. The
bizarre olympiad problem grew into the subject of {\em riemannian
security}. Namely, for a pair of (not necessarily distinct) points
$m_1,m_2$ in a riemannian manifold  $M$ let $\Ga(m_1,m_2)$ be the
set of geodesic segments joining these points. A set $B\subset
M\setminus\{m_1,m_2\}$ is {\em blocking} if every
$\ga\in\Ga(m_1,m_2)$ intersects $B$. The pair $m_1,m_2$ is {\em
secure} if there is a {\em finite blocking set} $B=B(m_1,m_2)$. A
manifold is secure if all pairs of points are secure. If there is
a uniform bound on the cardinalities of blocking sets, the
manifold is {\em uniformly secure} and the best possible bound is
the {\em blocking number}.

i) What closed riemannian manifolds are secure? ii) What plane
polygons are secure? It was the latter question that first got
into the literature \cite{HiSn98}. A polygon is {\em rational} if
its corners have $\pi$-rational angles. By \cite{HiSn98}, all
rational polygons are secure. The author studied the security of
{\em translation surfaces} \cite{GJ96,GJ00} and proved that the
regular $n$-gon is secure if and only if $n=3,4,6$ \cite{Gut05}.
Since all regular polygons are rational, this disproves the claim
in \cite{HiSn98}. The work \cite{Gut06} contains related results
on the security of rational polygons, but question ii) remains
wide open \cite{Gut12}. Question i) has been studied in
\cite{BanGut10,BG08,GeKu11,GeLi11,Gut09,GS06,HeKu12,Herr09,LS07,SchSou10}.
The following conjecture is widely accepted:
\begin{conj}   \label{secure_conj}
A closed riemannian manifold is secure if and only if it is flat.
\end{conj}
Flat manifolds are uniformly secure, and the blocking number
depends only on their dimension \cite{GS06,BanGut10,Gut05}. They
are also {\em midpoint secure}, i.e., the midpoints of connecting
geodesics yield a finite blocking set for any pair of points
\cite{GS06,BanGut10,Gut05}. Conjecture~~\ref{secure_conj} says
that flat manifolds are the only secure manifolds. This was
verified for several special cases: A manifold without conjugate
points is uniformly secure if and only if it is flat
\cite{BG08,LS07}; a compact locally symmetric space is secure if
and only if it is flat \cite{GS06}. The generic manifold is
insecure \cite{GeKu11,GeLi11,HeKu12}. Generic two-dimensional tori
are {\em totally insecure}, i.e., have no secure pairs of points
\cite{BanGut10}. Any riemannian metric has an arbitrarily close,
insecure metric in the same conformal class \cite{HeKu12}.
Riemannian surfaces of genus greater than one are totally insecure
\cite{BanGut10}.

\medskip

This paper adds evidence to the validity of
Conjecture~~\ref{secure_conj}, albeit indirectly. Integral curves
of a {\em spray} on a differentiable manifold play the role of
geodesics on a riemannian manifold \cite{Sh01}. In particular,
they yield the set of {\em connecting curves} for any pair of
points in $M$. This allows us to speak of (in)security for sprays
the same way we did for riemannian manifolds.

In this work we study this question for {\em Lie sprays} on
homogeneous spaces $M=G/\Ga$ where $G$ is a Lie group and
$\Ga\subset G$ is a lattice. Connection curves are the orbits of
one-parameter subgroups of $G$. To avoid confusion, we do not use
the term ``security'' in this setting. We speak of {\em finite
blocking} instead. The counterpart of ``secure'' in this context
is the term {\em blockable}. See Section~~\ref{setting}. The Lie
spray analog of Conjecture~~\ref{secure_conj} is as follows:
\begin{conj}   \label{block_conj}
Let $M=G/\Ga$ where $G$ is a connected Lie group and $\Ga\subset
G$ is a lattice. Then $M$ is blockable if and only if $G=\R^n$,
i.e., $M$ is a torus.
\end{conj}

Our main result, Theorem~~\ref{main_thm}, establishes
Conjecture~\ref{block_conj} for nilmanifolds. Minimal geodesics
proved to be a useful tool in riemannian security
\cite{Ba88,In86,Ba11}. The main tool in the present study is the
geometry of Lie groups \cite{Ma51,AusGrHa63,Star00,GorNev10}.
Section~~\ref{setting} recalls some properties of spaces $G/\Ga$.
Section~~\ref{heisen} characterizes blockable pairs of points in
nilmanifolds of the classical heisenberg group. In
Section~~\ref{general} we prove Conjecture~\ref{block_conj} for
nilmanifolds. In Section~~\ref{example} we characterize blockable
pairs of points in arbitrary heisenberg manifolds.
Section~~\ref{conj} reduces midpoint blocking in $\macn/\mczn$ to
a study of square roots of $\mczn$-cosets. We conclude with a
conjecture about midpoint blocking in $G/\Ga$ for simple
noncompact Lie groups.

\section{Connection blocking in homogeneous spaces} \label{setting}
We will study homogeneous spaces $M=G/\Ga$, where $G$ is a
connected Lie group, and $\Ga\subset G$ is a lattice.\footnote{Our
framework is valid for uniform and nonuniform lattices.} For $g\in
G,m\in M$ we denote by $g\cdot m$ the action of $G$. Let $\G$ be
the Lie algebra of $G$ and let $\exp:\G\to G$ be the exponential
map. For $m_1,m_2\in M$ let $C_{m_1,m_2}$ be the set of
parameterized curves $c(t)=\exp(tx)\cdot m,0\le t \le 1,$ such
that $c(0)=m_1,c(1)=m_2$. We say that $C_{m_1,m_2}$ is the
collection of {\em connecting curves} for the pair $m_1,m_2$. Let
$I\subset\R$ be any interval. If $c(t),t\in I,$ is a curve, we
denote by $c(I)\subset M$ the set $\{c(t):t\in I\}$. A {\em finite
set} $B\subset M\setminus\{m_1,m_2\}$ is a {\em blocking set} for
the pair $m_1,m_2$ if for any curve $c$ in $C_{m_1,m_2}$ we have
$c([0,1])\cap B\ne\emptyset$.

If a blocking set exists, we will say that {\em the pair $m_1,m_2$
is connection blockable}, often suppressing the adjective
`connection'. We will also say that $m_1$ is {\em blockable}
(resp. {\em not blockable})  {\em away} from $m_2$.

\begin{defin}   \label{blocking_def}
A homogeneous space $M=G/\Ga$ is {\em connection blockable} if
every pair of its points is blockable. If there exists at least
one non-blockable pair of points in $M$, then $M$ is
non-blockable.
\end{defin}

The analogy with riemannian security
\cite{Gut05,LS07,Herr09,BanGut12}\footnote{Some authors prefer to
use the terms `blocking' or `light blocking' in the riemannian
setting \cite{GeKu11,GeLi11,Herr09,LS07,SchSou10}.} suggests the
following:
\begin{defin}   \label{ramificat_def}
{\em 1.} A homogeneous space $M=G/\Ga$ is {\em uniformly
blockable} if there exists $N\in\N$ such that every pair of its
points can be blocked with a set $B$ of cardinality at most $N$.
The smallest such $N$ is the {\em blocking number} for $M$.

\noindent{\em 2.} A pair $m_1,m_2\in M$ is {\em midpoint
blockable} if the set $\{c(1/2):c\in C_{m_1,m_2}\}$ is finite. A
homogeneous space is midpoint blockable if all pairs of its points
are midpoint blockable.


\noindent{\em 3.} A homogeneous space is {\em totally
non-blockable} if no pair of its points is blockable.
\end{defin}
\begin{prop}  \label{prelim_prop}
Let $M=G/\Ga$ where $\Ga\subset G$ is a lattice, and let
$m_0=\Ga/\Ga\in M$. Then the following holds:

\noindent{\em 1.} The homogeneous space $M$ is  blockable (resp.
uniformly blockable, midpoint blockable) if and only if all pairs
$m_0,m$ are blockable (resp. uniformly blockable, midpoint
blockable). The space $M$ is totally non-blockable if and only if
no pair $m_0,m$ is blockable;

\medskip

\noindent{\em 2.} Let $\tGa\subset\Ga$ be lattices in $G$, let
$M=G/\Ga,\tM=G/\tGa$, and let $p:\tM\to M$ be the covering. Let
$m_1,m_2\in M$ and let $\tm_1,\tm_2\in\tM$ be such that
$m_1=p(\tm_1),m_2=p(\tm_2)$. If $B\subset M$ is a blocking set for
$m_1,m_2$ (resp. $\tB\subset\tM$ is a blocking set for
$\tm_1,\tm_2$) then $p^{-1}(B)$ (resp. $p(\tB)$) is a blocking set
for $\tm_1,\tm_2$ (resp. $m_1,m_2$).

\medskip

\noindent{\em 3.} Let $G',G''$ be connected Lie groups with
lattices $\Ga'\subset G',\Ga''\subset G''$, and let
$M'=G'/\Ga',M''=G''/\Ga''$. Set $G=G'\times G'',M=M'\times M''$.
Then a pair $(m_1',m_1''),(m_2',m_2'')\in M$ is connection
blockable if and only if both pairs $m_1',m_2'\in M'$ and
$m_1'',m_2''\in M''$ are connection blockable.
\begin{proof}
Claim 1 is immediate from the definitions. The proofs of claim 2
and claim 3 are analogous to the proofs of their counterparts for
riemannian security.  See Proposition 1 in \cite{GS06} for claim
2, and Lemma 5.1 and Proposition 5.2 in \cite{BG08} for claim 3.
\end{proof}
\end{prop}

Let $M_1,M_2$ be homogeneous spaces. We will use the following
terminology. Suppose that one of them is blockable (or not),
midpoint blockable (or not), totally non-blockable (or not), etc
if and only if the other one is. We will then say that both spaces
have {\em identical blocking properties}.

Recall that two subgroups $\Ga_1,\Ga_2\subset G$ are {\em
commensurable}, $\Ga_1\sim\Ga_2$, if there exists $g\in G$ such
that the group $\Ga_1\cap g\Ga_2 g^{-1}$ has finite index in both
$\Ga_1$ and $g\Ga_2 g^{-1}$. Commensurability yields an
equivalence relation in the set of lattices in $G$. We will use
the following immediate Corollary of
Proposition~~\ref{prelim_prop}.

\begin{cor}    \label{commensur_cor}
{\em 1.} If lattices $\Ga_1,\Ga_2\subset G$ are commensurable,
then the homogeneous spaces $M_i=G/\Ga_i:i=1,2$ have identical
blocking properties.

\noindent{\em 2.} Let $M_1=G_1/\Ga_1,M_2=G_2/\Ga_2$ be homogeneous
spaces. Let $M=M_1\times M_2=(G_1\times G_2)/(\Ga_1\times\Ga_2)$.
Then $M$ is blockable (resp. midpoint blockable, uniformly
blockable) if and only if both $M_1$ and $M_2$ are blockable
(resp. midpoint blockable, uniformly blockable).
\end{cor}

Let $\exp:\G\to G$ be the exponential map. For $\Ga\subset G$
denote by $p_{\Ga}:G\to G/\Ga$ the projection, and set
$\exp_{\Ga}=p_{\Ga}\circ\exp:\G\to G/\Ga$. We will say that a pair
$(G,\Ga)$ is of {\em exponential type} if the map $\exp_{\Ga}$ is
surjective. Let $M=G/\Ga$. For $m\in M$ set
$\Log(m)=\exp_{\Ga}^{-1}(m)$.

\begin{prop}  \label{criterion_prop}
Let $G$ be a Lie group, let $\Ga\subset G$ be a lattice such that
$(G,\Ga)$ is of exponential type, and let $M=G/\Ga$.

\noindent Then $m\in M$ is blockable away from $m_0$ if and only
if there is a map $x\mapsto t_x$ of $\Log(m)$ to $(0,1)$ such that
the set $\{\exp(t_xx):x\in\Log(m)\}$ is contained in a finite
union of $\Ga$-cosets.
\begin{proof}
Connecting curves are $c_x(t)=\exp(tx)\Ga/\Ga$ for some
$x\in\Log(m)$. Since $c(1)=m$, there is $\ga\in\Ga$ such that
$\exp(x)=g\ga$. Thus
\begin{equation}   \label{connect_eq}
c(t)=\exp(t\log(g\ga))\cdot m_0
\end{equation}
for some $\ga\in\Ga$, and every such curve is connecting $m_0$
with $m$.

Suppose $m$ is blockable away from $m_0$, and let $B\subset G/\Ga$
be a blocking set. Let $t_x\in(0,1)$ be such that $c_x(t_x)\in B$.
Set $A=\{\exp(t_xx):x\in\Log(m)\}\subset G$. Then
$(A\Ga/\Ga)\subset B$, hence finite. Thus $A$ is contained in a
finite union of $\Ga$-cosets.

On the other hand, if for any collection
$\{t_x\in(0,1):x\in\Log(m)\}$ the set
$A=\{\exp(t_xx):x\in\Log(m)\}$ is contained in a finite union of
$\Ga$-cosets, then $(A\Ga/\Ga)\subset M$ is a finite blocking set.
\end{proof}
\end{prop}

If $A\subset G$ any subset, we will say that
\begin{equation}   \label{sqrt_eq}
\Sqrt(A)=\{g\in G:g^2\in A\}.
\end{equation}
is the {\em square root} of $A$. We will say that a pair $(G,\Ga)$
is of {\em virtually exponential type} if there exists
$\tGa\sim\Ga$ such that $(G,\tGa)$ is of exponential type.

\begin{cor}    \label{criterion_cor}
Let $\Ga\subset G$ be a lattice such that $(G,\Ga)$ is of
virtually exponential type. Then:

\noindent{\em 1}. The homogeneous space $M=G/\Ga$ is midpoint
blockable if and only if the square root of any coset $g\Ga$ is
contained in a finite union of $\Ga$-cosets.

\noindent{\em 2}. Any point in $M$ is midpoint blockable away from
itself if and only if the square root of $\Ga$ is contained in a
finite union of $\Ga$-cosets.
\begin{proof}
By Corollary~~\ref{commensur_cor}, we can assume that $(G,\Ga)$ is
of exponential type. Set $t_x\equiv 1/2$ in
Proposition~~\ref{criterion_prop}.
\end{proof}
\end{cor}
\section{Connection blocking in three-dimensional heisenberg manifolds and some
other two-step nilmanifolds} \label{heisen}
For readers' convenience, we will recall basic facts about
connected, simply connected, real nilpotent Lie groups
\cite{Kiri62,AusGrHa63}. Let $G$ be as above. Its Lie algebra $\G$
has an ascending tower of ideals
$\{0\}\subset\G_1\subset\cdots\subset\G_{p-1}\subset\G_p=\G$ such
that $\G_i/\G_{i-1}$ is the center of $\G/\G_{i-1}$. We will say
that $\G$ (resp. $G$) is a {\em $p$-step nilpotent Lie algebra}
(resp. group). When $p=2$, the above decomposition becomes
$\CC\subset\G$ where $\CC$ is the center of $\G$ and $\G/\CC$ is
abelian.

The map $\exp:\G\to G$ is a diffeomorphism. Set $\log=\exp^{-1}$.
For $t\in\R$ we define the diffeomorphism $g\mapsto g^t$ of $G$ by
$g^t=\exp(t\log g)$. The haar measure in $G$ is both left and
right invariant. All lattices $\Ga\subset G$ are uniform
\cite{Ma51,AusGrHa63}. Referring to a measure on a nilmanifold
$M=G/\Ga$, we will always mean the unique invariant probability
measure $\mu$. By a measure on the set of pairs $m_1,m_2\in M$ we
will mean the measure $\mu\times\mu$.

\subsection{Blocking in the classical heisenberg manifold}  \label{classical_sub}
\hfill \break
The unique nonabelian nilpotent Lie algebra of $3$ dimensions is
$\hh=\R X+\R Y+\R Z$, where $[X,Y]=Z$ and $[X,Z]=[Y,Z]=0$. For
historical reasons, $\hh$ is usually called the {\em heisenberg
Lie algebra}. In the modern terminology, there is an infinite
sequence of heisenberg Lie algebras $\hh_n$, $n\ge 1$, where
$\hh_n$ is a two-step nilpotent Lie algebra of $2n+1$ dimensions.
The corresponding simply connected nilpotent groups $H_n,n\ge 1,$
are the {\em heisenberg groups}; the nilmanifolds $H_n/\Ga$ are
the {\em heisenberg manifolds}. In this subsection we study
blocking in a special heisenberg manifold.

We will denote $\hh_1$ and $H_1$ by $\hh$ and $H$ respectively. In
order to avoid confusion with the material in
section~~\ref{example}, we will speak of the classical heisenberg
group and the classical heisenberg manifold. It is standard to
represent $\hh$ and $H$ by $3\times 3$ matrices:
$$
\hh=\{\left[
\begin{array}{ccc}
0 & x & z \\
0 & 0 & y \\
0 & 0 & 0
\end{array}
\right]:x,y,z\in\R\},\ H=\{\left[
\begin{array}{ccc}
1 & x & z \\
0 & 1 & y \\
0 & 0 & 1
\end{array}
\right]:x,y,z\in\R\}.
$$
We will use the notation
$$
h(x,y,z)=\left[
\begin{array}{ccc}
1 & x & z \\
0 & 1 & y \\
0 & 0 & 1
\end{array}
\right].
$$

The {\em classical Heisenberg manifold} is $M=H/\Ga$ where
$\Ga=\{h(p,q,r):p,q,r\in\Z\}$. Using the unique decomposition
$h=h(a,b,c)h(p,q,r)$ where $0\le a,b,c<1,\,p,q,r\in\Z$, we
identify $M$ as a set with the unit cube $Q=[0,1)^3$. For
$(a,b,c)\in Q$ we denote by $m(a,b,c)\in M$ the corresponding
element. Then $m_0=m(0,0,0)$.

For $h=h(x,y,z)$ set
$\pi_x(h)=x\mod1,\pi_y(h)=y\mod1,\pi_z(h)=z\mod1$. Thus
$\pi_x,\pi_y,\pi_z:H\to\R/\Z$. We will denote by $\oplus$ the
addition in $\R/\Z$, i.e., $x\oplus y=x+y\mod1$. We will need a
criterion for a subset of $H$ to be contained in a finite inion of
$\Ga$-cosets.

\begin{lem}  \label{heisen_coset_lem}
Let $W\subset H$. Set $A=\pi_x(W),B=\pi_y(W),C=\pi_z(W)$.

\noindent Then $|W\Ga/\Ga|<\infty$ if and only if the sets $A,B$
are finite and $C\subset \cup_{i=1}^N\{c_i\oplus\Z a_i: a_i\in
A\}$.
\begin{proof}
The identity
\begin{equation}   \label{heisen_prod_eq}
h(a,b,c)h(p,q,r)=h(a+p,b+q,c+qa+r)
\end{equation}
implies that $W\subset h\Ga$ if and only if there exist
$a,b,c\in\R/\Z$ such that $A=\{a\},B=\{b\}$, and
$C\subset(c\oplus\Z a)$.
\end{proof}
\end{lem}

If $a,c\in\R/\Z$, we will refer to the set $c\oplus \Z
a\subset\R/\Z$ as a {\em rotation orbit}.

\begin{prop}  \label{heisen_gen_prop}
An element $m=m(a,b,c)\in M$ is blockable away from $m_0$ if and
only if $b\in\Q a+\Q$.
\begin{proof}
By a straightforward calculation
\begin{equation}   \label{heisen_power_eq}
\left(h(x,y,z)\right)^t=h(tx,ty,tz+\frac{t(t-1)}{2}xy).
\end{equation}
Let $(a,b,c)\in Q,(p,q,r)\in\Z^3$. Set
$h=h(a,b,c),\gamma=h(p,q,r)$, and let $0<t\le 1$.
Equation~~\eqref{heisen_power_eq} implies
\begin{equation}   \label{prod_power_eq}
\left(h\ga\right)^t=h(t(a+p),t(b+q),t(c+r+qa)+\frac{t(t-1)}{2}(a+p)(b+q)).
\end{equation}

By Proposition~~\ref{criterion_cor}, $m(a,b,c)$ is blockable away
from $m_0$ if and only if for each $(p,q,r)\in\Z^3$ there exist
$0<t_{pqr}<1$ such that the set
$$
W=\{\left(h(a,b,c)h(p,q,r)\right)^{t_{pqr}}:p,q,r\in\Z\}
$$
is contained in a finite union of $\Ga$-cosets. Set
$A=\pi_x(W),B=\pi_y(W),C=\pi_z(W)$. Then
$$
A=\cup_{p,q,r\in\Z}t_{pqr}a\oplus t_{pqr}p,\
B=\cup_{p,q,r\in\Z}t_{pqr}b\oplus t_{pqr}q
$$
and
$$
C=\{t_{pqr}c\oplus t_{pqr}r\oplus
t_{pqr}qa\oplus\frac{t_{pqr}(t_{pqr}-1)}{2}qa
$$
$$
\oplus\frac{t_{pqr}(t_{pqr}-1)}{2}pb
\oplus\frac{t_{pqr}(t_{pqr}-1)}{2}ab\oplus\frac{t_{pqr}(t_{pqr}-1)}{2}pq:p,q,r\in\Z\}.
$$
The sets $A$ and $B$ are finite if and only if
$T=\cup_{p,q,r\in\Z}\{t_{pqr}\}$ is a finite subset of $\Q$. Then
$C$ is contained in a finite union of rotation orbits by elements
in $A$ if and only if $b\in\Q a+\Q$. The claim now follows from
Lemma~~\ref{heisen_coset_lem}.
\end{proof}
\end{prop}
\begin{cor}   \label{heisen_thm}
Let $m_1=m(a_1,b_1,c_1),m_2=m(a_2,b_2,c_2)$ be arbitrary points in
the classical heisenberg  manifold. Then the pair $m_1,m_2$ is
blockable if and only if $b_1-b_2\in\Q(a_1-a_2)+\Q$.
\begin{proof}
Set $h=h(-a_1,-b_1,-c_1-a_1b_1)$. Then $h\cdot m_1=m_0,h\cdot
m_2=m(a_2-a_1,b_2-b_1,c_2-c_1-a_1(b_2+b_1))$. The claim follows
from the invariance of blockability under the group action and
Proposition~~\ref{heisen_gen_prop}.
\end{proof}
\end{cor}

We will now study connection blocking in nilmanifolds
$\tM=H/\tGa$, where $\tGa\subset H$ is an arbitrary lattice. If
the lattices $\Ga',\Ga''\subset H$ are isomorphic by an
automorphism of $H$, the nilmanifolds $H/\Ga',H/\Ga''$ have
identical blocking properties. Any lattice in $H$ is isomorphic by
an automorphism of $H$ to $\Ga(\de)=\{g(\de p,q,r):p,q,r\in\Z\}$
where $\de\in\N$ \cite{GoWi86}. Thus, it suffices to analyze
connection blocking in nilmanifolds $M_{\de}=H/\Ga(\de)$.

Set $Q_{\de}=[0,\de)\times[0,1)\times[0,1)$. The decomposition
$h=h(a,b,c)h(\de p,q,r)$, where $0\le a<\de,0\le
b,c<1,\,p,q,r\in\Z$, identifies $M_{\de}$ as a set with $Q_{\de}$.
For $(a,b,c)\in Q_{\de}$ we denote by $m^{(\de)}(a,b,c)\in
M_{\de}$ the corresponding element.

\begin{prop}  \label{heisen_more_gen_prop}
Let $m_1=m^{(\de)}(a_1,b_1,c_1),m_2=m^{(\de)}(a_2,b_2,c_2)$ be
arbitrary points in $M_{\de}$. Then the pair $m_1,m_2$ is
blockable if and only if $b_1-b_2\in\Q(a_1-a_2)+\Q$.
\begin{proof}
By Corollary~~\ref{heisen_thm}, the claim holds for $\de=1$. Let
$\de>1$. The inclusion $\Ga_{\de}\subset\Ga$ yields the
$\de$-to-$1$ covering $p_{\de}:M_{\de}\to M$ which, under the
identifications $M=Q,M_{\de}=Q_{\de}$, becomes
$p_{\de}(a,b,c)=(a\mod1,b,c)$.

Let $m_1,m_2\in M_{\de}$ be any pair. By
Proposition~~\ref{prelim_prop}, it is blockable if and only if the
pair $p_{\de}(m_1),p_{\de}(m_2)\in M$ is blockable. The claim now
follows from Corollary~~\ref{heisen_thm}.
\end{proof}
\end{prop}
\begin{thm}    \label{classical_heisen_cor}
Let $M$ be any three-dimensional heisenberg  manifold. Then the
following claims hold.

\noindent{\em 1.} A pair of points in $M$ is blockable if and only
if it is midpoint blockable.

\noindent{\em 2.} Every point in $M$ is blockable away from
itself.

\noindent{\em 3.}  The set of blockable pairs of points is a dense
countable union of closed submanifolds of positive codimension.
\begin{proof}
We assume without loss of generality that $M=M_{\de}$ for some
$\de\in\N$. Let $m_i=m^{(\de)}(a_i,b_i,c_i)$, $i=1,2$ be any pair
of points in $M$. By Proposition~~\ref{heisen_more_gen_prop}, the
pair $m_1,m_2$ is blockable if and only if
$b_1-b_2\in\Q(a_1-a_2)+\Q$. But then, by
Proposition~~\ref{prelim_prop} and
Proposition~~\ref{heisen_gen_prop}, it is midpoint blockable. This
proves claim 1. Claim 2 is immediate from
Proposition~~\ref{heisen_more_gen_prop}.

Let now $m_1=m^{(\de)}(a_1,b_1,c_1)\in M$. For any
$(a_2,c_2)\in[0,\de)\times[0,1)$ the set of numbers $b_2\in[0,1)$
such that $b_1-b_2\in\Q(a_1-a_2)+\Q$ is countably dense. By
Proposition~~\ref{heisen_more_gen_prop}, the set of elements
$m_2=m^{(\de)}(a_2,b_2,c_2)$ such that the pair $m_1,m_2$ is
blockable, is a dense countable union of two-dimensional
manifolds. Claim 3 follows.
\end{proof}
\end{thm}

Let $M=M_{\de}$ be a three-dimensional heisenberg manifold. The
identification $M_{\de}=Q_{\de}$ by $m=m^{(\de)}(a,b,c)$, sends
$\mu$ to the normalized lebesgue measure.

\begin{cor}    \label{measure_nul_cor}
Almost all pairs of points in a three-dimensional heisenberg
manifold are not blockable.
\begin{proof}
Let $(M\times M)_{\blo}\subset M\times M$ be the set of blockable
pairs. By the proof of Theorem~~\ref{classical_heisen_cor}, for
any $m_1\in M$ the set $\{m_2\in M:(m_1,m_2)\in(M\times
M)_{\blo}\}$ is a countable union of subsets of positive
codimension. Hence $\mu\left(\{m_2\in M:(m_1,m_2)\in(M\times
M)_{\blo}\}\right)=0$. The claim follows, by the Fubini theorem.
\end{proof}
\end{cor}
%


%
\subsection{Blocking in a family of two-step nilmanifolds}  \label{two_step_sub}
\hfill \break
Let $\G$ be a two-step nilpotent Lie algebra with the center $\CC$
such that $\dim(\G/\CC)=2$. Let $G$ be the connected, simply
connected Lie group with the Lie algebra $\G$.

\begin{prop}      \label{special_twostep_prop}
Let $C\subset G$ be the center, and let $d+1=\dim(C)$, where $d\ge
0$. Let $\Ga\subset G$ be a lattice. Set $N=G/\Ga$.

Then $N=M\times T^d$ where $M$ is a three-dimensional heisenberg
manifold and $T^d=\R^d/\Z^d$ is the $d$-dimensional torus. A pair
of points $n_1=(m_1,t_1),n_2=(m_2,t_2)$ in $N$ is blockable if and
only if the pair $m_1,m_2\in M$ is blockable.
\begin{proof}
If $d=0$ then $G=H_1$, and there is nothing to prove. Thus, we
assume from now on that $d\ge 1$. Let $Z\in\CC$ be the unique, up
to scalar multiple, element such that $Z=[X,Y]$ for some
$X,Y\in\G$. By \cite{Ma51}, we can choose elements $X,Y,Z$ so that
$\exp X,\exp Y,\exp Z\in\Ga$. Let $\V\subset\CC$ be the
$d$-dimensional subspace complementary to $\R Z$, spanned by
elements $v_1,\dots,v_d$ such that $\exp v_1,\dots,\exp
v_d\in\Ga$. Let $V=\exp\V\subset C$ and let $L=V\cap\Ga$. Then
$\exp:\G\to G$ induces the isomorphisms $V=\R^d,L=\Z^d$.

Set $\hh=\R X+\R Y+\R Z$ and $H=\exp\hh\subset G$. Then $H$ is the
classical heisenberg group, and $\Ga\cap H\subset H$ is a lattice.
Let $M=H/(\Ga\cap H)$ and $T^d=V/L$. The decompositions $G=H\times
V,\Ga=(\Ga\cap H)\times L$ yield the first claim.

By Corollary~~\ref{commensur_cor}, the pair $n_1,n_2$ is blockable
if and only if both pairs $m_1,m_2\in M$ and $t_1,t_2\in T^d$ are
blockable. Connection curves in $T^d$  are the geodesics for a
flat metric. Since a flat torus is secure \cite{GS06}, the claim
follows.
\end{proof}
\end{prop}
\begin{cor}    \label{special_twostep_cor}
Let $G$ be a two-step nilpotent Lie group with the center $C$
satisfying $\dim(G/C)=2$. Let $M$ be a $G$-nilmanifold. Then the
following claims hold.

\noindent{\em 1.} A pair of points in $M$ is blockable if and only
if it is midpoint blockable.

\noindent{\em 2.} Every point in $M$ is  blockable away from
itself.

\noindent{\em 3.} The set of blockable pairs of points is a dense
countable union of closed submanifolds of positive codimension.
\begin{proof}
Mimic the proof of Theorem~~\ref{classical_heisen_cor}, using
Proposition~~\ref{special_twostep_prop} instead of
Proposition~~\ref{heisen_more_gen_prop}.
\end{proof}
\end{cor}
\section{Connection blocking in arbitrary nilmanifolds}  \label{general}
Let $G$ be a connected, simply connected, nilpotent Lie group and
let $\Ga\subset G$ be a lattice. Let $M=G/\Ga$ be the
corresponding nilmanifold.

\begin{prop}      \label{general_prop}
If $G$ is not abelian, then there exist nonblockable pairs of
points in $M$.
\begin{proof}
Let $C\subset G$ be the center of $G$ and let $\CC\subset\G$ be
its Lie algebra. Then $\CC\subset\G$ is a proper inclusion, and
$\dim(\G)-\dim(\CC)\ge 2$ \cite{Kiri62}. There are $X,Y\in\G$ and
$Z\in\CC,Z\ne 0$ such that $[X,Y]=Z$ \cite{Kiri62}. Moreover, we
can choose these elements so that $\exp X,\exp Y,\exp Z\in\Ga$
\cite{Ma51}. Set $\G_1=\R X + \R Y + \CC$, and let $G_1\subset G$
be the corresponding subgroup. Then $G_1$ is a two-step nilpotent
Lie group with the center $C$ and $\dim(G_1/C)=2$. The group
$\Ga_1=\Ga\cap G$ is a lattice in $G_1$. Let $M_1=G_1/\Ga_1$ be
the corresponding nilmanifold.

Let $m_0\in M$ be the base point. Then $G_1\cdot m_0=M_1\subset
M$. If a pair $m_1,m_2\in M_1$ is not blockable in $M_1$, then it
is not blockable in $M$. The claim now follows from
Corollary~~\ref{special_twostep_cor}.
\end{proof}
\end{prop}

\medskip

\begin{thm}     \label{main_thm}
Let $M$ be a nilmanifold of $n$ dimensions. Then the following
statements are equivalent:

\noindent {\em 1.} The manifold is connection blockable;

\noindent {\em 2.} The manifold is midpoint blockable;

\noindent {\em 3.} We have $\pi_1(M)=\Z^n$;

\noindent {\em 4.} It is a topological torus;

\noindent {\em 5.} It is uniformly blockable; the blocking number
depends only on its dimension.
\begin{proof}
Let $M=G/\Ga$. By Proposition~~\ref{general_prop}, if $G$ is not
abelian, then $M$ is not blockable. If $G$ is abelian, then
$M=\R^n/\Z^n$; connecting curves are the geodesics in any flat
riemannian metric. Since flat tori are secure \cite{GS06}, $M$ is
connection blockable. Thus, $M$ is connection blockable if and
only if $M=\R^n/\Z^n$.

Malcev \cite{Ma51} proved that compact nilmanifolds are isomorphic
if and only if they are homemorphic if and only if they have the
same fundamental group. Thus, statements $1$, $3$, and $4$ are
equivalent. Flat manifolds, in particular, flat tori are midpoint
secure \cite{Gut05,BanGut10,BanGut12}. The canonical parameter for
connecting curves is proportional to the arc length parameter.
Thus, tori are midpoint blockable, proving the equivalence of
statements $1$ and $2$. The implication $1\to 5$ is a consequence
of the observation that the blocking number for flat tori of $n$
dimensions is $2^n$ \cite{Gut05} and the Bieberbach theorem
\cite{BanGut10}.
\end{proof}
\end{thm}
\begin{cor}    \label{mid_block_cor}
Each point in a nilmanifold $M^n$ is blockable away from itself.
Moreover, the blocking is uniform, and the optimal bound depends
only on $n$.
\begin{proof}
Let $M=G/\Ga$ be any homogeneous manifold. By
Corollary~~\ref{criterion_cor}, either all points in $M$ are
blockable away from themselves or no point in $M$ is blockable
away from itself. The former happens if and only if
$|\Sqrt(\Ga)\Ga/\Ga|<\infty$. For lattices in nilpotent Lie groups
this property holds \cite{Ma51}. Moreover, there exist $c_n\in\N$
such that for any lattice $\Ga$ in a nilpotent Lie group $G$ of
$n$ dimensions, $|\Sqrt(\Ga)\Ga/\Ga|<c_n$ \cite{Ma51}. The second
claim follows, by Proposition~~\ref{criterion_prop}.
\end{proof}
\end{cor}
\section{Blocking in arbitrary Heisenberg manifolds}  \label{example}
Let $n\ge 1$. For $\vx,\vy\in\R^n,z\in\R$ set
$$
h_n(\vx,\vy,z)=\left[
\begin{array}{ccc}
1 & \vx & z \\
0 & \id_n & \vy \\
0 & 0 & 1
\end{array}
\right],
$$
where $\vx$ (resp. $\vy$) is the row (resp. column) vector. The
group $H_n=\{h_n(\vx,\vy,z):\vx,\vy\in\R^n,z\in\R\}$ is the
$(2n+1)$-dimensional heisenberg group. Heisenberg manifolds are
the nilmanifolds $H_n/\Ga$ where $\Ga\subset H_n$ is a lattice.
For $\vde=(\de_1,\dots,\de_n)\in\Z_+^n$ let
$\vde\Z^n=\{(\de_1k_1,\dots,\de_nk_n):\vk\in\Z^n\}$. The group
$\Ga_n(\vde)=\{h_n(\vx,\vy,z):\vx\in\vde\Z^n,\vy\in\Z^n,z\in\Z\}\subset
H_n$ is a lattice. Set $M_n(\vde)=H_n/\Ga_n(\vde)$. We will first
study blocking in nilmanifolds $M_n=M_n(\vde)$ for $\vde=\vec{1}$.
Using the unique decomposition $h=h(\va,\vb,c)h(\vvp,\vq,r)$ where
$ \va,\vb\in[0,1)^n,c\in[0,1),\,\vvp,\vq\in\Z^n,r\in\Z$, we
identify $M_n$ as a set with the $(2n+1)$-dimensional cube
$Q_{2n+1}=[0,1)^n\times[0,1)^n\times[0,1)$. For $(\va,\vb,c)\in
Q_{2n+1}$ we denote by $m(\va,\vb,c)\in M_n$ the corresponding
element. Then $m_0=m(\vec{0},\vec{0},0)$.

For $h=h(\vx,\vy,z)$ set
$\pi_x(h)=\vx\mod\vec{1},\pi_y(h)=\vy\mod\vec{1},\pi_z(h)=z\mod1$.
Thus $\pi_x,\pi_y:H\to\R^n/\Z^n,\pi_z:H\to\R/\Z$. By $\oplus$ we
will denote the addition in $\R^k/\Z^k$ for any $k$.

\begin{lem}  \label{multi_heisen_lem}
Let $W\subset H_n$. Set $A=\pi_x(W),B=\pi_y(W),C=\pi_z(W)$.

\noindent Then $|W\Ga/\Ga|<\infty$ if and only if the sets $A,B$
are finite and
$C\subset\cup_{i=1}^N\{c_i\oplus\Z\langle\vq_i,\va_i\rangle:\va_i\in
A,\vq_i\in\Z^n\}$.
\end{lem}

The proof of Lemma~~\ref{multi_heisen_lem} is the obvious
modification of the argument in Lemma~~\ref{heisen_coset_lem}, and
we leave it to the reader. By Lemma~~\ref{multi_heisen_lem}, a set
$W\subset G$ satisfies $|W\Ga/\Ga|<\infty$ if and only if the sets
$A,B$ are finite, and the set $C$ is contained in a finite union
of orbits of rotation by $\langle\vq,\va\rangle,\va\in A$. We
denote by $L(n,\R),L(n,\Q),L(n,\Z)$ the sets of $n\times n$
matrices with entries in $\R,\Q,\Z$ respectively.

\begin{prop}  \label{multi_heisen_prop}
A pair $m_1=m(\va_1,\vb_1,c_1),m_2=m(\va_2,\vb_2,c_2)$ in $M$ is
blockable if and only if there exist a matrix $L\in L(n,\Q)$ and a
vector $\vl\in\Q^n$ such that $\vb_1-\vb_2=L(\va_1-\va_2)+\vl$.
\begin{proof}
Mimic the proofs of Corollary~~\ref{heisen_thm}, replacing
Lemma~~\ref{heisen_coset_lem} by Lemma~~\ref{multi_heisen_lem}.
\end{proof}
\end{prop}

We will now study connection blocking in nilmanifolds $M_n(\vde)$
for arbitrary $\vde$. Set
$[\vec{0},\vde)=\prod_{i=1}^n[0,\de_i)\subset\R_+^n$ and
$Q_{2n+1}(\vde)=[\vec{0},\vde)\times[0,1)^n\times[0,1)$. The
unique decomposition $h=h(\va,\vb,c)\ga$ where $(\va,\vb,c)\in
Q_{2n+1}(\vde),\ga\in\Ga_n(\vde)$, identifies $M_n(\vde)$ as a set
with $Q_{2n+1}(\vde)$. For $(\va,\vb,c)\in Q_{2n+1}(\vde)$ we
denote by $m^{(\vde)}(\va,\vb,c)\in M_n(\vde)$ the corresponding
element. Then $m_0^{(\vde)}=m^{(\vde)}(\vec{0},\vec{0},0)$.

\begin{prop}  \label{multiheisen_more_gen_prop}
Let
$m_1=m^{(\vde)}(\va_1,\vb_1,c_1),m_2=m^{(\vde)}(\va_2,\vb_2,c_2)$
be arbitrary points in $M_n(\vde)$. Then the pair $m_1,m_2$ is
blockable if and only if there exist a matrix $L\in L(n,\Q)$ and a
vector $\vl\in\Q^n$ such that $\vb_1-\vb_2=L(\va_1-\va_2)+\vl$.
\begin{proof}
Set $|\vde|=\de_1\cdots\de_n$. Let $\Ga_n$ be the standard integer
lattice in $H_n$. The inclusion $\Ga_n(\vde)\subset\Ga_n$ yields
the $|\vde|$-to-$1$ covering $M_n(\vde)\to M_n$. Now mimic the
proof of Proposition~~\ref{heisen_more_gen_prop}, replacing
Corollary~~\ref{heisen_thm} by
Proposition~~\ref{multi_heisen_prop}.
\end{proof}
\end{prop}

The following theorem summarizes the properties of connection
blocking in heisenberg manifolds.

\begin{thm}    \label{gener_heisen_cor}
Let $M$ be any heisenberg manifold. Then the following claims
hold.

\noindent{\em 1.} A pair of points in $M$ is blockable if and only
if it is midpoint blockable.

\noindent{\em 2.} Every point in $M$ is blockable away from
itself.

\noindent{\em 3.}  The set of blockable pairs of points is a dense
countable union of closed submanifolds of positive codimension.

\noindent{\em 4.} Almost all pairs of points are not blockable.
\begin{proof}
By \cite{GoWi86}, any Lattice $\Ga\subset H_n$ is isomorphic to
$\Ga_n(\vde)$ by an automorphism of $H_n$. Hence, it suffices to
prove the claims for the nilmanifolds $M_n(\vde)$. The arguments
are the multidimensional versions of those used to prove
Theorem~~\ref{classical_heisen_cor} and
Corollary~~\ref{measure_nul_cor}, with
Proposition~~\ref{heisen_more_gen_prop} replaced by
Proposition~~\ref{multiheisen_more_gen_prop}. Details are left to
the reader.
\end{proof}
\end{thm}
\section{Blocking in semi-simple homogeneous manifolds: Examples and conjectures}  \label{conj}
We will illustrate connection blocking in homogeneous spaces
$G/\Ga$ which are not nilmanifolds with an example and will
formulate a conjecture. We will need the following Lemma. The
proof is straightforward, and we leave it to the reader.

\begin{lem}     \label{subgroup_lem}
Let $G$ be a Lie group, and let $\Ga\subset G$ be a lattice. Let
$H\subset G$ be a closed subgroup such that $\Ga\cap H$ is a
lattice in $H$. Let $X=G/\Ga,Y=H/(\Ga\cap H)$ be the homogeneous
spaces, and let $Y\subset X$ be the natural inclusion.

\noindent {\em 1}. If $Y$ is not connection blockable (resp. not
midpoint blockable, etc) then $X$ is not connection blockable
(resp. not midpoint blockable, etc). {\em 2}. If $Y$ contains a
point which is not blockable (resp. midpoint blockable) away from
itself, then no point in $X$ is blockable (resp. midpoint
blockable) away from itself.
\end{lem}

The following Lemma will be used in the proof of
Proposition~~\ref{slnR_exa}.

\begin{lem}     \label{sl2_root_lem}
Let $\left[\begin{array}{cc}
                  a  & b  \\
                  c  & d   \\
\end{array}\right],\ \left[\begin{array}{cc}
                  x  & y  \\
                  z  & w   \\
\end{array}\right]\in\mac$. If $a+1,d+1,a+d+2,c\ne 0$, then
\begin{equation}   \label{square_eq1}
\left[\begin{array}{cc}
                  x  & y  \\
                  z  & w   \\
\end{array}\right]^2 = \left[\begin{array}{cc}
                  a  & b  \\
                  c  & d   \\
\end{array}\right]
\end{equation}
if and only if
\begin{equation}   \label{square_eq2}
z^2=\frac{c^2}{a+d+2},\ x=\frac{a+1}{c}z,\ y=\frac{b}{c}z,\
w=\frac{d+1}{c}z.
\end{equation}
\begin{proof}
Equation~~\eqref{square_eq1} is equivalent to
$$
\left[\begin{array}{cc}
                  x  & y  \\
                  z  & w   \\
\end{array}\right] = \left[\begin{array}{cc}
                  w  & -y  \\
                  -z  & x   \\
\end{array}\right]
\left[\begin{array}{cc}
                  a  & b  \\
                  c  & d   \\
\end{array}\right].
$$
Solve for $x,y,w$ in terms of $z$ and use $xw-yz=1$.
\end{proof}
\end{lem}

We will now study connection blocking in the homogeneous space
$\macn/\mczn$.

\begin{prop}    \label{slnR_exa}
The space $\macn/\mczn$ contains pairs of points which are not
midpoint blockable. In particular, no point in $m\in\macn/\mczn$
is midpoint blockable away from itself.
\begin{proof}
Set $M_n=\macn/\mczn$. Although the exponential map for $\macn$ is
not surjective, the pair $(\macn,\mczn)$ is of exponential type.
Let $m_1,m_2\in M_n$, let $g\in\macn$ satisfy $m_2=g\cdot m_1$. By
Proposition~~\ref{criterion_prop} and
Corollary~~\ref{criterion_cor}, the pair $m_1,m_2$ is midpoint
blockable if and only if $\Sqrt(g\cdot\mczn)$ is contained in a
finite union of $\mczn$-cosets.

Let $g=\left[\begin{array}{cc}
                  a  & b  \\
                  c  & d   \\
\end{array}\right],\ X=\left[\begin{array}{cc}
                  x  & y  \\
                  z  & w   \\
\end{array}\right]\in\mac,\ \left[\begin{array}{cc}
                  p  & q  \\
                  r  & s   \\
\end{array}\right]\in\mcz$. By Lemma~~\ref{sl2_root_lem},
if $X\in\Sqrt(g\mcz)$ then
$z^2=(pc+rd)^2(pa+rb+qc+sd+2)^{-1}$. Let $K$ be the smallest field
containing $a,b,c,d$. By the above, the entries of matrices $X$ in
$\Sqrt(g\mcz)$ contain the square roots of infinitely many
$\Z$-independent elements in $K$. Hence, the $\Z$-module generated
by these matrix entries  has infinite $\Z$-rank. On the other
hand, the $\Z$-module generated by the entries of matrices in a
finite union of $\mcz$-cosets has finite $\Z$-rank. Therefore,
$\Sqrt(g\cdot\mcz)$ is not contained in a finite union of
$\mcz$-cosets.

\medskip

For $n>2,\,1\le i\le n-1$ let $G_i\subset\macn$ be the group
$\mac$ imbedded in $\macn$ via the rows and columns $i,i+1$. Then
$G_i\cap\mczn=\mcz$, and hence $G_i\mczn/\mczn=\mac/\mcz$. Set
$M_n^{(i)}=G_i\mczn/\mczn\subset M_n$. By
Lemma~~\ref{subgroup_lem}, no pair $m_1,m_2$ in $M_n^{(i)}$ is
midpoint blockable, yielding the former part of the claim.

By Proposition~~\ref{criterion_prop} and
Corollary~~\ref{criterion_cor}, elements in $M_n$ are midpoint
blockable away from themselves if and only if $\Sqrt(\mczn)$ is
contained in a finite union of $\mczn$-cosets. Since the identity
element belongs to all $G_i$, the set $\Sqrt(\mczn)$ is not
contained in a finite union of $\mczn$-cosets, yielding the claim.
\end{proof}
\end{prop}

The preceding argument shows that no pair $m_1,m_2$ in $\mac/\mcz$
is midpoint blockable. We have not shown this for $\macn/\mczn$ if
$n>2$. However, the above proof and related considerations suggest
the following.

\begin{conj}     \label{mid_block_conj}
Let $G$ be a  simple, connected, noncompact Lie group, and let
$\Ga\subset G$ be a lattice. Then no pair in $G/\Ga$ is midpoint
blockable.
\end{conj}

\medskip

\noindent Acknowledgements. The work was partially supported by
the MNiSzW grant N N201 384834 and the NCN Grant
DEC-2011/03/B/ST1/00407.

\end {document}